\documentclass[a4paper, 10pt, journal, onecolumn]{ieeeconf}  

\IEEEoverridecommandlockouts                              

\usepackage{cite}
\usepackage{mathtools}
\usepackage{ifamath} 
\usepackage{subcaption}
\let\labelindent\relax 
\usepackage[inline]{enumitem}
\usepackage{xcolor,soul}
\usepackage{pgfplots}
\usepgfplotslibrary{fillbetween}
\usepackage{tikz}
\usepackage{tikzscale}
\usepackage[]{units}
\usetikzlibrary{calc,external,patterns}
\definecolor{blue}{rgb}{0.25,0.25,0.75}
\colorlet{someblue}{blue!70}
\colorlet{lightblue}{blue!55}
\colorlet{verylightblue}{blue!30}
\definecolor{greyblue}{rgb}{0.36,0.36,0.5}
\colorlet{lightgreyblue}{greyblue!55}
\colorlet{red}{red!90!black}
\colorlet{lightred}{red!65}
\colorlet{verylightred}{red!40}
\colorlet{green}{green!80!black}
\colorlet{lightgreen}{green!65}
\colorlet{verylightgreen}{green!20}
\colorlet{lightgray}{gray!65}

\usepackage[%
  hidelinks=true,
  pdftitle={Robust Control Policies given Formal Specifications in Uncertain Environments},
  pdfauthor={Damian Frick,  Tony A. Wood, Gian Ulli, Maryam Kamgarpour},
  pdfcreator={pdfLaTeX},
  pdfsubject={}]{hyperref}

\newcommand{\apSet}{{\rm AP}}

\title{\LARGE \bf
Robust Control Policies given Formal Specifications in Uncertain Environments$^\star$}

\author{Damian Frick, Tony~A. Wood, Gian Ulli, Maryam Kamgarpour
\thanks{$^\star$ This manuscript is the preprint of a paper published in IEEE Control System Letters and is subject to IEEE Control Systems Society copyright. The IEEE Control Systems Society maintains the sole rights of distribution or publication of the work in all forms and media. The copy of record is available at \href{https://doi.org/10.1109/LCSYS.2017.2700333}{10.1109/LCSYS.2017.2700333}.}
\thanks{The authors are with the Automatic Control Laboratory, Department of Electrical Engineering and Information Technology, ETH Zurich, 8092 Zurich, Switzerland. E-mail addresses: {\tt \{dafrick,woodt,mkamgar\}@control.ee.ethz.ch}~for~\{D.~Frick, T.~A. Wood, M.~Kamgarpour\} and {\tt ug@student.ethz.ch} for G.~Ulli.}%
\thanks{This research was supported by the Swiss National Science Foundation (Synergia) No. 141836 and European Union ERC Starting Grant, CONENE.}
}

\begin{document}
\maketitle

\thispagestyle{empty}
\pagestyle{empty}

\begin{abstract}
We consider robust control synthesis for linear systems with complex specifications that are affected by uncertain disturbances. This work is motivated by autonomous systems interacting with partially known, time-varying environments. Given a specification in bounded linear temporal logic, we propose a method to synthesize control policies that fulfill the specification for all considered disturbances, where the disturbances affect the atomic propositions of the specification. Our approach relies on introducing affine disturbance feedback policies and casting the problem as a bilinear optimization problem. We introduce an inner approximation of the constraint set leading to a mixed-integer quadratic program that can be solved using general-purpose solvers. The framework is demonstrated on a numerical case study with applications to autonomous driving.
\end{abstract}

\section{Introduction}

Increased autonomy in applications ranging from transportation to energy systems necessitates the synthesis of controllers that perform safely in the presence of uncertainties. Often, such control laws need to satisfy complex specifications. To move beyond single objective and point-to-point motion planning towards complex specifications, formal methods, such as linear temporal logic (LTL) \cite{pnueli1977}, can be used. LTL combines propositional logic with temporal operators to enable the expression of logical statements in time.
Once a specification is formulated, standard model checking tools can be used to synthesize hybrid controllers \cite{fainekos2005b,tabuada2006,belta2007} based on a finite-state abstraction which bisimulates the original continuous system.
For discrete-time linear or mixed-logical dynamical systems \cite{bemporad1999}, the problem of finding trajectories that satisfy LTL specifications can be posed as a mixed-integer linear or quadratic program (MILP/MIQP) \cite{karaman2008,wolff2014} and can be applied in a receding horizon fashion \cite{gol2015}. This approach naturally allows for the consideration of continuous states and inputs. It benefits from the computational advances in mixed-integer optimization algorithms, rather than constructing a possibly large discrete abstraction.

Dealing with uncertainties is an area of active research in the context of control synthesis with formal specifications. The source of such uncertainties can be broadly classified into two categories: \emph{internal}, meaning uncertainty in the system dynamics or the system model, and \emph{external}, meaning uncertainty affecting the environment. Such external uncertainties can correspond to uncertainty in target locations or uncertain behavior of the environment. Furthermore, the objective when dealing with uncertainty is usually two-fold \begin{enumerate*}[label=(\roman*)]
	\item to improve controller performance by incorporating knowledge about the uncertainty into the controller synthesis, e.g., taking into account known probability distributions of, or bounds on, the disturbance, and
	\item to robustify the controller against the disturbances, avoiding violation of the specifications for any disturbance realization.
\end{enumerate*}

Stochastic models have been used to cope with probabilistic uncertainty affecting the dynamical system, whereas probabilistic specification languages have been introduced to address external uncertainty. In \cite{kamgarpour2013} stochastic hybrid systems and a subset of LTL specifications are considered. Maximizing the probability of satisfying the specification is cast as a stochastic reachability problem. This is generalized to include probabilistic uncertainties in the location of goal and obstacle sets \cite{kamgarpour2017}. In \cite{ding2011,wolff2012} general LTL specifications for a finite-state Markov decision process (MDP) are considered.
A dynamic programming approach is proposed to synthesize controllers that satisfy the specification.
Finally in \cite{lahijanian2012} MDPs are used with a probabilistic specification language.
These approaches aim to achieve satisfaction of the specification in probability and do not consider other performance criteria.
This is in contrast to \cite{aoude2013}, where a minimum-time objective is pursued, while ensuring that uncertain obstacles are avoided with a given probability. However, \cite{aoude2013} does not address general LTL specifications.

In a \emph{robust} setting, past work has addressed uncertain environments.  
A powerful concept for temporal logic planning in uncertain environments is \emph{reactive planning/synthesis}. In this framework, generalized reactivity(1) specifications \cite{piterman2006} are used that capture both the task specifications and the allowed uncertain behavior of the environment. In \cite{kress-gazit2009} standard tools for LTL controller synthesis are used to generate controllers that are robust to uncertain environment behavior. In addition, \cite{wongpiromsarn2012} addresses the receding horizon case. However, these methods cannot easily capture dynamic disturbances.
Alternatively, in \cite{raman2015} signal temporal logic, a more advanced specification language that captures robustness, is used. It quantifies the degree of satisfaction of a specification. Trajectories that maximize robustness can be generated via repeated solution of MILPs.

\subsection{Contribution}
In this work, we focus on LTL specifications in uncertain or time-varying environments. Our motivation stems from the presence of uncertainties in the obstacles or goal sets. Hence, we introduce a framework in which the atomic propositions in LTL are themselves affected by uncertainty. 
In contrast to methods presented in the literature, we synthesize a \emph{robust} control \emph{policy}.
We propose a novel approach based on the use of affine disturbance feedback policies to cope with uncertainties. The use of feedback improves performance compared to open-loop control policies by taking into account measurements of past disturbances in real time. While this approach is well-known in the model predictive control literature, its application to formal method control synthesis, to the best of our knowledge, has not been explored before. To deal with the computational complexity of the resulting robust MILP/MIQP, we propose an inner approximation and illustrate its performance via a case study.

\section{Problem Formulation} \label{sec:problem}

We consider discrete-time linear systems
\begin{equation} \label{eqn:sys}
	x_{k+1} = Ax_{k} + B u_{k}\,,
\end{equation}
where $x_k \in \reals{n_x}$ is the system state at time $k$ and $u_k \in \reals{n_u}$ is the control input applied between time $k$ and $k+1$.

Note that results of this work extend to systems affected by additive disturbances and can further be extended to discrete-time hybrid dynamics described via mixed-integer constraints in the framework of mixed logical dynamical systems \cite{bemporad1999}.

\subsection{Uncertain Temporal Logic Specifications}

In safety critical planning problems it is often desirable to impose strict specifications for the allowed trajectories. Specifications can include statements such as \emph{reach-avoid}, reaching a goal set while avoiding obstacles, or \emph{coverage}, visiting a collection of regions.
Linear temporal logic allows the rigorous description of such specifications.
For system~\eqref{eqn:sys}, we define a finite trajectory, or \emph{run}, of length $L$ starting at $x_j$, as
\begin{equation*}
\mathbf{x}_j^L := \begin{bmatrix}x_j^\transp & x_{j+1}^\transp & \ldots & x_{j+L}^\transp\end{bmatrix}^\transp\,,
\end{equation*}
a sequence of states $x_k$ such that for each $k=j,\ldots,j+L-1$ there exists an input $u_k$ such that $x_{k+1} = Ax_{k} + Bu_{k}$.

We consider LTL specifications that are affected by an uncertain \emph{disturbance} vector $w \in \reals{n_w}$. Given a specification $\varphi$, length $L$ and index $j$, we want to find a control input sequence $\mathbf{u}_j^L := \begin{bmatrix}u_j^\transp & \cdots & u_{j+L-1}^\transp\end{bmatrix}^\transp$ such that the run $\mathbf{x}_j^L$ satisfies the specification $\varphi$ for all disturbance sequence realizations $\mathbf{w}_j^L := \begin{bmatrix}w_j^\transp & \ldots & w_{j+L}^\transp\end{bmatrix}^\transp$ contained in a bounded polyhedron $\setbf{W}_j^L \subseteq \reals{(L+1)n_w}$, i.e., 
\begin{equation*}
	(\mathbf{x}_j^L,\mathbf{w}_j^L) \satisfies \varphi\quad \forall \mathbf{w}_j^L \in \setbf{W}_j^L\,.
\end{equation*}

For simplicity we consider LTL formulae in positive normal form \cite{baier2008}. To avoid issues with unbounded effects of the disturbances, we furthermore use bounded LTL formulae, without loops \cite[Definition~2.1]{biere2006}, a subset of the usual LTL semantics.
A formula in LTL is a combination of \emph{atomic propositions} $p$ taken from a finite set $\apSet := \{p_1, \ldots , p_m\}$, propositional logic operators $\lnot$~(\emph{not}), $\land$~(\emph{and}) and $\lor$~(\emph{or}), and temporal operators $\lnext$~(\emph{next}), $\until$~(\emph{until}) and $\release$~(\emph{release}). More formally, we define LTL formulae via the grammar
\begin{equation*}
	p \sep{} \lnot p \sep{} \phi \land \psi \sep{} \phi \lor \psi \sep{} \lnext \phi \sep{} \phi \until \psi \sep{} \phi \release \psi\,,
\end{equation*}
where $\phi, \psi$ are LTL formulae. Atomic propositions take values in $\{\true,\false\}$. In the context of this work, the disturbance $w$ enters the description of the atomic propositions $p_i \in \apSet$, i.e., each $p_i$ is associated with a polyhedral set
\begin{equation*}
	\set{P}_i := \{(x,w) \in \reals{n_x+n_w} \sep{} P_i^x x \leq  P_i^w w + \rho_i \}\,,
\end{equation*}
defined over the state-disturbance space.

Given a sequence of disturbance realizations $\mathbf{w}_j^L$, a run $\mathbf{x}_j^L$ satisfies an atomic proposition $p_i$ if the \emph{augmented state} $\mathbf{z}_j^L:=(z_j,\ldots,z_{j+L})$, with $z_{{\hat{\jmath}}} := (x_{\hat{\jmath}},w_{\hat{\jmath}})$, satisfies $z_j \in \set{P}_i$. The satisfaction of the formula $p_i$ is denoted by $\mathbf{z}_j^L \satisfies p_i$. The propositional operators are defined as
\begin{subequations} \label{eqn:LTLSemantics}
\begin{align}
	&\mathbf{z}_j^L \satisfies \lnot p_i &&\text{ iff } z_j \not\in \set{P}_i\,,\label{eqn:LTLSemanticsNot}\\
	&\mathbf{z}_j^L \satisfies \phi \land \psi &&\text{ iff } \mathbf{z}_j^L \satisfies \phi \text{ and } \mathbf{z}_j^L \satisfies \psi\,,\label{eqn:LTLSemanticsAnd}\\
	&\mathbf{z}_j^L \satisfies \phi \lor \psi &&\text{ iff } \mathbf{z}_j^L \satisfies \phi \text{ or } \mathbf{z}_j^L \satisfies \psi\,,\label{eqn:LTLSemanticsOr}\\
	\intertext{and the temporal operators are defined as}
	&\mathbf{z}_j^L \satisfies \lnext \phi &&\text{ iff } \mathbf{z}_{j+1}^{L-1} \satisfies \phi\,, \label{eqn:LTLSemanticsNext}\\
	&\mathbf{z}_j^L \satisfies \phi \until \psi &&\text{ iff } \exists \hat{\jmath} \in \{0, \ldots, L-1\}\suchthat \mathbf{z}_{j+\hat{\jmath}}^{L-\hat{\jmath}} \satisfies \psi \text{ and } \forall i \in \{0,\ldots,\hat{\jmath}-1\}\,: \mathbf{z}_{j+i}^{L-i} \satisfies \phi\,,\label{eqn:LTLSemanticsUntil}\\
	&\mathbf{z}_j^L \satisfies \phi \release \psi &&\text{ iff } \forall \hat{\jmath} \in \{0,\ldots,L-1\}:\mathbf{z}_{j+\hat{\jmath}}^{L-\hat{\jmath}} \satisfies \psi \text{ or } \exists i\in\{0,\ldots,\hat{\jmath}-1\} \suchthat \mathbf{z}_{j+i}^{L-i} \satisfies \phi\,.\label{eqn:LTLSemanticsRelease}
\end{align}
\end{subequations}
We introduce the additional temporal operators $\eventually \phi := \true \until \phi$ (\emph{eventually}) and $\always \phi :=  \false \release \phi$ (\emph{always}).

\subsection{Robust Policy Synthesis}
Given a fixed \emph{planning horizon} $N$ and \emph{initial state} $x_0$, we define the state trajectory $\mathbf{x} := \mathbf{x}_0^N$, disturbance sequence $\mathbf{w} := \mathbf{w}_0^N$ and corresponding input $\mathbf{u} := \mathbf{u}_0^N$, with
\begin{equation*}
	\setbf{W} := \setbf{W}_0^N := \{\mathbf{w} \in \reals{(N+1)n_w} \sep{} \mathbf{W}\mathbf{w} \leq \mathbf{v} \}\,,
\end{equation*}
a closed and \emph{bounded} polyhedron, with $\mathbf{W} \in \reals{n_v \times (N+1)n_w}$ and $\mathbf{v} \in \reals{n_v}$.

Employing causal disturbance feedback policies allows us to synthesize a control law that can react to past disturbances in real-time based on measured data. Such robust feedback policies can be generated, even though only bounds on the disturbances are known during control synthesis.
\begin{problem}
Given an LTL specification $\varphi$, find a sequence of causal disturbance feedback policies
\begin{equation*}
	u_0(w_0),\ldots,u_{N-1}(w_0,\ldots,w_{N-1})\,,
\end{equation*}
such that for all realizations of the uncertainty $\mathbf{w} \in \setbf{W}$:
we minimize an objective function and satisfy
\begin{enumerate*}[label=(\roman*)]
	\item input constraints $u_k \in \set{U} \subseteq \reals{n_u}$,
	\item state constraints $x_{k+1} \in \set{X} \subseteq \reals{n_x}$, and
	\item the specification $(\mathbf{x},\mathbf{w}) \satisfies \varphi$
\end{enumerate*}.
\end{problem}

\section{Solution Approach} \label{sec:solution}
Searching for general feedback policies is intractable. Hence, we focus on linear feedback.
We define a causal, affine disturbance feedback policy with parameters $\mathbf{H} \in \reals{Nn_u\times Nn_{w}}$ and $\mathbf{h} \in \reals{Nn_u}$:
\begin{equation} \label{eqn:affpolicy}
	\mathbf{u} = 
	\underbrace{\begin{bmatrix}H_{1,1} & & & 0\\
	                           \vdots & \ddots & & \vdots\\
	                           H_{N,1} & \cdots & H_{N,N} & 0 \end{bmatrix}}_{\begin{bsmallmatrix}\mathbf{H} & 0\end{bsmallmatrix}}\mathbf{w}
	             + \underbrace{\begin{bmatrix}h_1 \\ \vdots \\ h_N\end{bmatrix}}_{\mathbf{h}}\,.
\end{equation}
Such affine policies are used in robust control \cite{bemporad1998, zhang2017} to improve performance compared to open-loop policies. By setting $\mathbf{H} = 0$ we have an open-loop policy, which is similar to the approach in \cite{raman2015}. There, worst-case disturbance sequences are computed via MILPs and then the trajectory is robustified against those sampled sequences. In our approach, using linear programming duality, only one mixed-integer program (MIP) needs to be solved to obtain guarantees for all disturbances.

Substituting policy \eqref{eqn:affpolicy} into the discrete-time system equations \eqref{eqn:sys}, we express the trajectory $\mathbf{x}$ as a function of the initial state $x_0$, the disturbance sequence $\mathbf{w}$, and the parameters $\mathbf{H}$ and $\mathbf{h}$ of our policy:
\begin{align*}
	\mathbf{x} &= \underbrace{\begin{bmatrix}A^0 \\ A^1 \\ \vdots \\ A^N\end{bmatrix}}_{\mathbf{A}}x_0 + \underbrace{\begin{bmatrix}0 \\ B \\ AB & B \\ \vdots & & \ddots \\ A^{N-1}B & \cdots & AB & B \end{bmatrix}}_{\mathbf{B}}\mathbf{u}\\
	&= \mathbf{A}x_0 + \begin{bmatrix}\mathbf{B}\mathbf{H} & 0\end{bmatrix}\mathbf{w} + \mathbf{B}\mathbf{h}\,.
\end{align*}
We consider the following robust policy synthesis problem:
\begin{problem}[Robust policy synthesis] \label{prob:aff}
\setlength{\abovedisplayskip}{0.3em} \setlength{\abovedisplayshortskip}{0.3em}
\begin{align}
	\min_{\mathbf{H},\mathbf{h}} \: & J(\mathbf{H},\mathbf{h})\nonumber\\
	\suchthat & \begin{rcases}\begin{bmatrix}\mathbf{H} & 0\end{bmatrix}\mathbf{w} + \mathbf{h} \in \setbf{U}\,,\\
	\mathbf{A}\theta + \begin{bmatrix}\mathbf{B}\mathbf{H} & 0\end{bmatrix}\mathbf{w} + \mathbf{B}\mathbf{h} \in \setbf{X}\,,\\
	(\mathbf{A}\theta + \begin{bmatrix}\mathbf{B}\mathbf{H} & 0\end{bmatrix}\mathbf{w} + \mathbf{B}\mathbf{h},\mathbf{w}) \satisfies \varphi
	\end{rcases} \forall \mathbf{w} \in \setbf{W}\,,\label{eqn:specconst}
\end{align}
where, $J : \reals{Nn_u\times N n_{w}} \times \reals{Nn_u} \rightarrow \reals{}$ is the convex quadratic objective and $\theta \in \reals{n_x}$ is the parametric initial state. Furthermore, $\setbf{U} := \set{U} \times \ldots \times \set{U} \subseteq \reals{Nn_u}$ and $\setbf{X} := \set{X} \times \ldots \times \set{X} \subseteq \reals{Nn_x}$ are the ``stacked'' input and state constraints, with $\set{U}$ and $\set{X}$ compact, convex polyhedra.
\end{problem}

\begin{remark}
	The case where the disturbance has known linear dynamics $v_{k+1} = A^w v_k +B^w w_k$ with an uncertain input $w_k$, naturally fits into the presented framework. Disturbance state feedback $u_k = K_{k+1} v_k + \kappa_{k+1}$ can be equivalently posed as disturbance feedback by appropriate choice of $\setbf{W}$ and restrictions on the structure of $\mathbf{H}$.
\end{remark}

\subsection{Conversion to Robust Mixed-Integer Program} \label{sec:mipconversion}
It is known that specification constraints of the form $\mathbf{z} \satisfies \varphi$ can be transformed into linear mixed-integer inequalities \cite{wolff2014}. This is achieved by introducing auxiliary variables, some of which are restricted to be binary. We transform the constraint
\begin{equation*}
	(\mathbf{A}\theta + \begin{bmatrix}\mathbf{B} \mathbf{H} & 0\end{bmatrix}\mathbf{w} + \mathbf{B}\mathbf{h},\mathbf{w}) \satisfies \varphi\,,
\end{equation*}
of \refprob{prob:aff} into a set of mixed-integer inequalities:
\begin{align} \label{eqn:miineq}
	f^\varphi(\mathbf{H},\mathbf{h},\theta,\mathbf{w},\delta) :=& \left(F^x\begin{bmatrix}\mathbf{B}\mathbf{H} & 0\end{bmatrix} + F^w\right)\mathbf{w} + F^x\mathbf{B}\mathbf{h} + F^x\mathbf{A}\theta + F^\delta \delta + f \leq 0\,,
\end{align}
where the auxiliary vector $\delta \in \Delta := \reals{n_c} \times \binaries{n_b}$ consists of $n_c$ continuous and $n_b$ binary variables.
The matrices $F^x$, $F^w$, $F^\delta$ and vector $f$ are of appropriate dimensions.
Notice, that \eqref{eqn:miineq} is linear in $\theta$, $\mathbf{h}$ and $\delta$, and bilinear in $\mathbf{H}$ and $\mathbf{w}$.
We augment the specification constraint \eqref{eqn:miineq} with the input and state constraints of \refprob{prob:aff}. This yields the following set of mixed-integer inequalities:
\begin{align} \label{eqn:ltl+ineq}
	g^\varphi(\mathbf{H},\mathbf{h},\theta,\mathbf{w},\delta) :=& \left(G^H\begin{bmatrix}\mathbf{B}\mathbf{H} & 0\end{bmatrix} + G^w\right)\mathbf{w} + G^H\mathbf{B}\mathbf{h} + G^\theta\mathbf{A}\theta + G^\delta \delta + g \leq 0\,,
\end{align}
where $g^\varphi(\mathbf{H},\mathbf{h},\theta,\mathbf{w},\delta)$ consists of $m^\varphi$ constraints and the matrices $G^H$, $G^w$, $G^\theta$ and $G^\delta$, as well as the vector $g$ have appropriate dimensions.
The number $n_b$ of binary variables does not change.
Substituting the set of mixed-integer inequalities \eqref{eqn:ltl+ineq} into \refprob{prob:aff} yields an optimization problem with linear mixed-integer constraints that need to be satisfied robustly, i.e., for all $\mathbf{w} \in \setbf{W}$. To make this problem tractable, we will reduce the robust constraint to a set of mixed-integer constraints.

\subsection{Reduction to Mixed-Integer Program}

The set of robustly admissible feedback gains $\mathbf{H}$ and $\mathbf{h}$, parametrized by the initial state $\theta$, is given as follows:
\begin{align} \label{eqn:miineq+}
	\set{C}^\varphi :=& \big\{ (\mathbf{H},\mathbf{h},\theta) \sep{\big} \eqref{eqn:specconst} \text{ holds for } \mathbf{H},\mathbf{h} \text{ given } \theta \, \big\}\nonumber\\
	=& \big\{ (\mathbf{H},\mathbf{h},\theta) \sep{} \forall\mathbf{w} \in \setbf{W}\:\exists \delta \in \Delta \suchthat g^\varphi(\mathbf{H},\mathbf{h},\theta,\mathbf{w},\delta) \leq 0 \, \big\}\nonumber\\
	=& \big\{ (\mathbf{H},\mathbf{h},\theta) \sep{} \hspace*{-0.4em} \max_{\mathbf{w}\in\setbf{W}} \min_{\delta \in \Delta} \max_{i} g^\varphi_i(\mathbf{H},\mathbf{h},\theta,\mathbf{w},\delta) \leq 0 \, \big\},
\end{align}
where $g^\varphi_i$ denotes the $i$-th component of $g^\varphi$. The solution of $\min_{\delta \in \Delta} \max_{i} g^\varphi_i(\mathbf{H},\mathbf{h},\theta,\mathbf{w},\delta)$ is a piecewise affine function in $\mathbf{H}\mathbf{w}$, $\mathbf{h}$ and $\theta$ \cite{ramirez2001}. The challenge in solving \eqref{eqn:miineq+} is due to the max-min structure and the bilinear dependence on $\mathbf{H}$ and $\mathbf{w}$. In particular, the maximization over $\mathbf{w}$ leads to non-linear constraints \cite{khalilpour2014}. Using linear programming duality, $\set{C}^\varphi$ can be represented by constraints that are linear in $\mathbf{h}$ and $\theta$, but contain bilinear terms involving~$\mathbf{H}$. This results in a representation of \refprob{prob:aff} as a bilinear program, which can be solved via spatial branch-and-bound \cite{androulakis1995}. For the case of $\mathbf{H} = 0$, the set $\set{C}^\varphi$ can be represented by linear mixed-integer inequalities, leading to a mixed-integer formulation which can be solved using general-purpose MIQP solvers.

To avoid the difficulty of solving a bilinear optimization problem for the general case $\mathbf{H}\neq 0$, we propose a simple inner approximation. This approximation produces a mixed-integer program and usually helps to preserve sparsity in the optimization problem. We consider the set
\begin{align} \label{eqn:innerapprox}
	\overline{\set{C}}^\varphi := \Big\{ (\mathbf{H},\mathbf{h},\theta) \sep{\Big} &\exists \delta \in \Delta\,,\:\greekbf{\lambda}_i \in \posreals{n_v} \text{ for } i=1,\ldots,m^\varphi \suchthat \text{for all } i = 1,\ldots,m^\varphi:\nonumber\\[-0.2em]
	& \quad \mathbf{v}^\transp \greekbf{\lambda}_i + G^H_{i\bcdot}\mathbf{B}\mathbf{h} + G_{i\bcdot}^\theta\mathbf{A}\theta + G_{i\bcdot}^\delta \delta + g_i \leq 0\,,\nonumber\\[-0.2em]
	& \quad \mathbf{W}^\transp \greekbf{\lambda}_i = \left(G^H_{i\bcdot}\begin{bmatrix}\mathbf{B}\mathbf{H} & 0\end{bmatrix} + G_{i\bcdot}^w\right)^\transp \Big\}\,,
\end{align}
where $G_{i\bcdot}$ denotes the $i$-th row of matrix $G$.
\begin{lemma} \label{lem:inner}
	The constraint set $\overline{\set{C}}^\varphi$ is an inner approximation of the robust specification constraint set $\set{C}^\varphi$, i.e., any feasible policy that satisfies \eqref{eqn:innerapprox} also satisfies \eqref{eqn:miineq+}.\\
\begin{proof}
First, we exchange the maximization over $\mathbf{w}$ with the minimization over $\delta$ in \eqref{eqn:miineq+}. This leads to the inner approximation
\begin{align*}
	\overline{\set{C}}^\varphi := \big\{ (\mathbf{H},\mathbf{h},\theta) \sep{\big} &\exists \delta \in \Delta \suchthat \hspace*{-0.8em} \max_{\substack{\mathbf{w}\in\setbf{W}\\i=1,\ldots,m^\varphi}} \hspace*{-0.5em}  g^\varphi_i(\mathbf{H},\mathbf{h},\theta,\mathbf{w},\delta) \leq 0 \big\} \subseteq \set{C}^\varphi\nonumber\,.
\end{align*}
For a given $i \in 1,\ldots,m^\varphi$, we use a standard robust optimization technique~\cite[p.~472]{bertsimas2011}, applying linear programming duality to replace $\max_{\mathbf{w}\in\setbf{W}} g^\varphi_i(\mathbf{H},\mathbf{h},\theta,\mathbf{w},\delta)$ with its dual
\begin{align*}
	\min_{\greekbf{\lambda}_i \in \posreals{n_v}} & \mathbf{v}^\transp \greekbf{\lambda}_i + G^H_{i\bcdot}\mathbf{B}\mathbf{h} + G_{i\bcdot}^\theta\mathbf{A}\theta + G_{i\bcdot}^\delta \delta + g_i\\[-0.8em]
	\suchthat & \mathbf{W}^\transp \greekbf{\lambda}_i = \left(G^H_{i\bcdot}\begin{bmatrix}\mathbf{B}\mathbf{H} & 0\end{bmatrix} + G_{i\bcdot}^w\right)^\transp\,.
\end{align*}
Dropping the minimization over $\greekbf{\lambda}_i$ gives an upper bound. Collecting these upper bounds for $i = 1,\ldots,m^\varphi$ yields \eqref{eqn:innerapprox}.
\end{proof}
\end{lemma}

Using the inner approximation $\overbar{\set{C}}^\varphi$ in \refprob{prob:aff} results in \refprob{prob:affmip}, a mixed-integer program that can be solved using general-purpose MIQP solvers.
\begin{problem}[MIQP] \label{prob:affmip}
\setlength{\abovedisplayskip}{0.3em} \setlength{\abovedisplayshortskip}{0.3em}
\begin{align*}
	\min_{\mathbf{H},\mathbf{h},\delta,\greekbf{\lambda}} \: & J(\mathbf{H},\mathbf{h})\\
	\suchthat & \mathbf{\overbar{V}}^\transp \greekbf{\lambda} + G^H\mathbf{B}\mathbf{h} + G^\theta\mathbf{A}\theta + G^\delta \delta + g \leq 0\,,\\
	& \overbar{\mathbf{W}}^\transp \greekbf{\lambda} = (I_{m^\varphi} \otimes \begin{bmatrix}\mathbf{H} & 0\end{bmatrix}^\transp) \mathbf{G}^\transp + {\mathbf{g}^\transp}\,,\\
	& \delta \in \Delta\,, \greekbf{\lambda} \in \posreals{m^\varphi n_v}\,,
\end{align*}
where $\otimes$ denotes the Kronecker product and $\greekbf{\lambda} := \begin{bmatrix}\greekbf{\lambda}_1^\transp & \ldots & \greekbf{\lambda}_{m^\varphi}^\transp\end{bmatrix}^\transp$, $\mathbf{\overbar{V}} := \diag(\mathbf{v}, \ldots, \mathbf{v})$,
$\overbar{\mathbf{W}} := \diag(\mathbf{W},\ldots,\mathbf{W})$,
$\mathbf{G} := \begin{bmatrix}G^H_{1\bcdot} \mathbf{B} & \cdots & G^H_{m^\varphi\bcdot} \mathbf{B} \end{bmatrix}$ and 
$\mathbf{g} := \begin{bmatrix}G^w_{1\bcdot} & \cdots & G^w_{m^\varphi\bcdot} \end{bmatrix}$.
\end{problem}
\begin{theorem} \label{thm:mip}
Any solution $\mathbf{H}^\star$, $\mathbf{h}^\star$ of \refprob{prob:affmip} is a feasible, possibly suboptimal, solution of \refprob{prob:aff}.\\
\begin{proof}
The proof follows directly from \reflem{lem:inner}.
\end{proof}
\end{theorem}
From \refthm{thm:mip} it follows that the policies obtained by solving \refprob{prob:affmip} robustly satisfy the specification for all disturbance realizations.
\begin{remark}
\refprob{prob:affmip} has $O(m^\varphi n_v + N^2 n_u n_w + n_c)$ continuous and $n_b$ binary decision variables, as well as $O(m^\varphi n_w)$ constraints. Mixed-integer solvers that can be used for instances of \refprob{prob:affmip} have a worst-case exponential complexity in the number of binary variables, $n_b$, whereas the complexity is polynomial, usually cubic, in the number of continuous variables and constraints. As remarked in \refsec{sec:mipconversion}, $n_b$ depends linearly on $N$. Solving \refprob{prob:affmip} therefore requires time that is in the worst-case exponential in the planning horizon $N$.
\end{remark}

\subsection{On the role of feedback and the approximation scheme}
The purpose of this section is to illustrate two points. First, in \refexa{exa:ex1}, we illustrate a simple case in which the inner approximation is not tight. In \refexa{exa:ex2}, we illustrate a case in which no feasible open-loop policy exists but a feasible disturbance feedback policy can be found.
\begin{figure*}
	\centering
	\begin{minipage}[b]{0.73\columnwidth}
	\begin{subfigure}[b]{0.31\linewidth}
	\centering
	\begin{tikzpicture}[scale=1,cap=round]
	\tikzstyle{dotted}=[dash pattern=on 0.5pt off 2pt]
	\begin{axis}[
		every outer x axis line/.append style={white!20!black},
		every x tick label/.append style={font=\color{white!20!black}},
		xmin=-1.5,
		xmax=1.5,
		xtick={-1,0,1},
		ymin=-1,
		ymax=1,
		ytick={-2,-1,0,1},
		axis lines = middle,
  		enlargelimits = true,
		xlabel={$x{-}w$},
		x label style={at={(axis description cs:1.2,0.25)}},
		ylabel={$e^\varphi(x,w)$},
		y label style={at={(axis description cs:0.5,1.2)}},
		axis on top,
		width = 1.8in, height = 1.4in]
		\addplot[name path=mfl,domain=-1.2:-0.5,color=blue,line width=1.5pt] {-x-1}; \addplot[name path=mfl,domain=-1.5:-1.2,color=blue,line width=1.5pt, dotted] {-x-1};
		\addplot[name path=mfl,domain=-0.5:0,color=blue,line width=1.5pt] {x};
		\addplot[name path=mfl,domain=0:0.5,color=blue,line width=1.5pt] {-x};
		\addplot[name path=mfl,domain=0.5:1.2,color=blue,line width=1.5pt] {x-1}; \addplot[name path=mfl,domain=1.2:1.5,color=blue,line width=1.5pt, dotted] {x-1};
    \end{axis}
	\end{tikzpicture}
	\caption{Piecewise affine constraint function $e^\varphi(H,h,\theta,w)$ for $x = \theta + Hw + h$ and $w$.\\\\}\label{fig:ex1:exact:pwa}
	\end{subfigure}
	\hspace*{0.02\linewidth}
	\begin{subfigure}[b]{0.31\linewidth}
	\centering
	\begin{tikzpicture}[scale=1,cap=round]
	\tikzstyle{dotted}=[dash pattern=on 0.5pt off 2pt]
	\begin{axis}[
		every outer x axis line/.append style={white!20!black},
		every x tick label/.append style={font=\color{white!20!black}},
		xmin=-1.5,
		xmax=1.5,
		xtick={-1,0,1},
		ymin=-1,
		ymax=1,
		ytick={-2,-1,0,1},
		axis lines = middle,
  		enlargelimits = true,
		xlabel={$h{+}\theta$},
		x label style={at={(axis description cs:1.2,0.22)}},
		ylabel={$\overbar{e}^\varphi(1,h,\theta,\delta)$},
		y label style={at={(axis description cs:0.5,1.2)}},
		axis on top,
		width = 1.8in, height = 1.4in]
		\addplot[name path=mfl,domain=-1.2:-0.5,color=blue,line width=1.5pt] {-x-1}; \addplot[name path=mfl,domain=-1.5:-1.2,color=blue,line width=1.5pt, dotted] {-x-1};
		\addplot[name path=mfl,domain=-0.5:0.2,color=blue,line width=1.5pt] {x}; \addplot[name path=mfl,domain=0.2:0.5,color=blue,line width=1.5pt,dotted] {x};
		\addplot[name path=mfl,domain=-0.2:0.5,color=lightred,line width=1.5pt, dashed] {-x}; \addplot[name path=mfl,domain=-0.5:-0.2,color=lightred,line width=1.5pt,dotted] {-x};
		\addplot[name path=mfl,domain=0.5:1.2,color=lightred,line width=1.5pt, dashed] {x-1}; \addplot[name path=mfl,domain=1.2:1.5,color=lightred,line width=1.5pt, dotted] {x-1};
    \end{axis}
	\end{tikzpicture}
	\caption{Mixed-integer constraint function $\overbar{e}^\varphi(H,h,\theta,\delta)$ for $H=1$. Cases $\delta = 1$ (\textbf{\color{blue}solid}), $\delta = 0$ (\textbf{\color{lightred}dashed}) illustrated for values of $h$ and $\theta$}\label{fig:ex1:inner:pwa}
	\end{subfigure}
	\hspace*{0.02\linewidth}
	\begin{subfigure}[b]{0.31\linewidth}
	\centering
	\begin{tikzpicture}[scale=1,cap=round]
	\begin{axis}[
		every outer x axis line/.append style={white!20!black},
		every x tick label/.append style={font=\color{white!20!black}},
		xmin=-1,
		xmax=1,
		xtick={-1,0,1},
		xticklabels={${-}1{-}\theta$,${-}\theta$,$1{-}\theta$},
		ymin=0,
		ymax=2.2,
		ytick={1, 2},
		axis lines = middle,
  		enlargelimits = true,
		xlabel={$h$},
		x label style={at={(axis description cs:1.1,0.1)}},
		ylabel={$H$},
		y label style={at={(axis description cs:0.5,1.1)}},
		axis on top,
		after end axis/.code={\path (axis cs:0,0) node [anchor=north,yshift=-0.2em] {{\color{white!20!black}${-}\theta$}};},
		width = 1.8in]
		\addplot[name path=mfl,domain=-1:0,color=blue,line width=1.5pt] {-x};
		\addplot[name path=mfu,domain=-1:0,color=blue,line width=1.5pt] {x+2};
		\addplot[name path=pfu,domain=0:1,color=blue,line width=1.5pt] {2-x};
		\addplot[name path=pfl,domain=0:1,color=blue,line width=1.5pt] {x};
		\addplot[fill=verylightblue] fill between [of=mfl and mfu,soft clip={domain=-1:0}];
		\addplot[fill=verylightblue] fill between [of=pfl and pfu,soft clip={domain=0:1}];
		\addplot[name path=mfl,domain=-1:-0.5,color=lightred,line width=1.5pt,dashed] {-x};
		\addplot[name path=mfu,domain=-1:-0.5,color=lightred,line width=1.5pt,dashed] {x+2};
		\addplot[name path=pfu,domain=-0.5:0,color=lightred,line width=1.5pt,dashed] {1-x};
		\addplot[name path=pfl,domain=-0.5:0,color=lightred,line width=1.5pt,dashed] {1+x};
		\addplot[fill=verylightred] fill between [of=mfl and mfu,soft clip={domain=-1:-0.5}];
		\addplot[fill=verylightred] fill between [of=pfl and pfu,soft clip={domain=-0.5:0}];
		\addplot[name path=mfl,domain=0:0.5,color=lightred,line width=1.5pt,dashed] {1-x};
		\addplot[name path=mfu,domain=0:0.5,color=lightred,line width=1.5pt,dashed] {x+1};
		\addplot[name path=pfu,domain=0.5:1,color=lightred,line width=1.5pt,dashed] {2-x};
		\addplot[name path=pfl,domain=0.5:1,color=lightred,line width=1.5pt,dashed] {x};
		\addplot[fill=verylightred] fill between [of=mfl and mfu,soft clip={domain=0:0.5}];
		\addplot[fill=verylightred] fill between [of=pfl and pfu,soft clip={domain=0.5:1}];
    \end{axis}
	\end{tikzpicture}
	\caption{The set $\set{C}^\varphi$ of all feasible affine disturbance feedback policies ({\color{blue}\textbf{solid}}) and its inner approximation $\overbar{\set{C}}^\varphi$ ({\color{lightred}\textbf{dashed}}). Both are parametrized by $\theta$.}\label{fig:ex1:feas}
	\end{subfigure}
	\caption{Illustrating the difference between the exact robust specification constraint formulation $\set{C}^\varphi$ and the mixed-integer inner approximation $\overbar{\set{C}}^\varphi$ for \refexa{exa:ex1}.} \label{fig:ex1}
	\end{minipage}
	\hspace*{0.02\columnwidth}
	\begin{minipage}[b]{0.23\columnwidth}
	\centering
	\begin{tikzpicture}[scale=1,cap=round]
	\begin{axis}[
		every outer x axis line/.append style={white!20!black},
		every x tick label/.append style={font=\color{white!20!black}},
		xmin=-4,
		xmax=4,
		xtick={-4,-2,0,2,4},
		xticklabels={$-4$,$ $,$0$,$ $,$4$},
		ymin=-4,
		ymax=4,
		axis y line=none,
		axis x line=middle,
  		enlargelimits = true,
		axis on top,
		width = 1.8in]
		\filldraw[blue,fill=verylightblue,line width=1.5pt] (axis cs:-2,-2) rectangle (axis cs:2,2);
		\addplot[color=red,mark=x,line width = 1.5pt,forget plot]
		table[row sep=crcr]{0    0\\};
		\node[right] at (axis cs:0.1,0.5) {{\color{red}$x_0$}};
		\draw[->,someblue,line width=1.5pt] (axis cs:2,0) -- (axis cs:3,0) node[near end, above,someblue] {$w_0$};
		\draw[->,someblue,line width=1.5pt] (axis cs:-2,0) -- (axis cs:-3,0);
		\node at (axis cs:0,-3) {$k=0$};
    \end{axis}
	\end{tikzpicture}\\
	\vspace*{-0.5cm}
	\begin{tikzpicture}[scale=1,cap=round]
	\begin{axis}[
		every outer x axis line/.append style={white!20!black},
		every x tick label/.append style={font=\color{white!20!black}},
		xmin=-4,
		xmax=4,
		xtick={-4,-2,0,2,4},
		xticklabels={$-4$,$ $,$0$,$ $,$4$},
		ymin=-4,
		ymax=4,
		axis y line=none,
		axis x line=middle,
  		enlargelimits = true,
		x label style={at={(axis description cs:1.1,0.1)}},
		y label style={at={(axis description cs:0.5,1.1)}},
		axis on top,
		width = 1.8in]
		\filldraw[blue,fill=verylightblue,line width=1.5pt] (axis cs:-1,-2) rectangle (axis cs:3,2);
		\addplot[color=red,mark=x,line width = 1.5pt,forget plot]
		table[row sep=crcr]{1    0\\};
		\draw[->,black,line width=1.5pt] (axis cs:0,0) -- (axis cs:1,0) node[near start, above,black] {$u_0$};
		\node[right] at (axis cs:1.1,0.5) {{\color{red}$x_1$}};
		\draw[->,someblue,line width=1.5pt] (axis cs:3,0) -- (axis cs:4,0) node[near end, above,someblue] {$w_1$};
		\draw[->,someblue,line width=1.5pt] (axis cs:-1,0) -- (axis cs:-2,0);
		\node at (axis cs:0,-3) {$k=1$};
    \end{axis}
	\end{tikzpicture}
	\caption{Illustration of robustly feasible feedback policy for \refexa{exa:ex2}.} \label{fig:ex2}
	\end{minipage}
\end{figure*}

\begin{example} \label{exa:ex1}
	We consider a planning horizon of $N=1$ and a simple dynamical system $x = \theta + u$ with initial state $\theta$. A scalar disturbance $w \in \set{W}:= [-1,1]$ affects the specification $\varphi := p \lor q$, where the polyhedra associated with $p$ and $q$ are
	\begin{align*}
		\set{P} &:= \{x \in \reals{} \sep{} x \in [-1,0]+w \} \text{ and } \set{Q} := \{x \in \reals{} \sep{} x \in [0,1]+w \}\,.
	\end{align*}
	Given the affine disturbance feedback policy $u : = Hw + h$, we can construct the piecewise affine constraint function
	\begin{equation*}
		e^\varphi(H,h,\theta,w) := \min_{\delta \in \Delta} \max_{i} g^\varphi_i(H,h,\theta,w,\delta)\,,
	\end{equation*}
	illustrated in \reffig{fig:ex1:exact:pwa}. We see that $e^\varphi(H,h,\theta,w) \leq 0$, i.e., the constraint is feasible, for all $x$ and $w$ such that $-1 \leq x-w \leq 1$. This yields a description of the set $\set{C}^\varphi$ with $H$ and $h$, parametrized by $\theta$ and illustrated ({\color{blue}\textbf{solid}}) in \reffig{fig:ex1:feas}. In this example $\set{C}^\varphi$ is a convex polyhedron. However, in general it may be neither polyhedral nor convex.
	We also consider the inner approximation $\overbar{\set{C}}^\varphi$ with the corresponding mixed-integer constraint function
	\begin{equation*}
		\overbar{e}^\varphi(H,h,\theta,\delta) := \max_{w \in \set{W}} \max_{i} g^\varphi_i(H,h,\theta,w,\delta)\,,
	\end{equation*}
	illustrated in \reffig{fig:ex1:inner:pwa} for $H=1$. For $H \neq 1$ the illustration in \reffig{fig:ex1:inner:pwa} needs to be shifted upwards by $|H-1|$, reducing the available choices for $h$. The set $\overbar{\set{C}}^\varphi$ of $H$ and $h$, parametrized by $\theta$, is depicted ({\color{red}\textbf{dashed}}) in \reffig{fig:ex1:feas}. 
	We see that $\overbar{\set{C}}^\varphi \subset \set{C}^\varphi$. Furthermore, when $H=0$, the inner approximation is empty, while the exact solution has a unique feasible policy: $h = -\theta$.
\end{example}

\begin{example}  \label{exa:ex2}
	We consider simple integrator dynamics $x_{k+1} = x_k + u_k$ with initial state $x_0 = 0$. We define a safe set $\set{P}_{\rm safe} := \{x \in \reals{2} \sep{} x_1 \in [-2,2]+w\,, x_2 \in [-2,2] \}$ with associated atomic proposition $p_{\rm safe}$. The disturbance $w_k \in \reals{}$ is in $[-1,1]$ for all time steps $k$. The initial state is illustrated in \reffig{fig:ex2}. The objective is to find an input sequence $u_0,u_1,\ldots$ that satisfies the specification $\varphi := \always p_{\rm safe}$ for all disturbance realizations and all time steps $k$. Clearly this is not possible, because the safe set can move both either left or right, i.e., no open-loop policy that robustly satisfies this specification exists. However, it is easy to see, that the affine disturbance feedback policy $u_k = w_k$ is robustly feasible.
\end{example}

%
%

\section{Case Study} \label{sec:study}

We consider a motion planning task on a two-lane highway illustrated in \reffig{fig:cs:track}. Two cars with known initial positions are cruising on the upper lane at uncertain velocities in the range $[\unitfrac[93.4]{km}{h},\unitfrac[106.6]{km}{h}]$. The controlled car is driving on the lower lane and has an initial velocity of $\unitfrac[110]{km}{h}$. It is modeled as a 2-dimensional double-integrator affected by a constant drag term:
\begin{equation} \label{eqn:doubleint}
	\begin{bmatrix}\dot{x}_1 \\ \dot{x}_2\end{bmatrix} = \begin{bmatrix}x_3 \\ x_4\end{bmatrix}\,, \text{ and } \begin{bmatrix}\dot{x}_3 \\ \dot{x}_4\end{bmatrix} = \begin{bmatrix}u_1 \\ u_2\end{bmatrix} - c_d\,,
\end{equation}
where the state $x \in \reals{4}$ contains the position $(x_1,x_2)$ of the car and the longitudinal and lateral velocities $(x_3,x_4)$. The accelerations $(u_1,u_2)$ are the control inputs and are limited to the set $\set{U}$: $u_1 \in [\unitfrac[-4]{m}{s^2},\unitfrac[2]{m}{s^2}]$ and $u_2 \in [\unitfrac[-3]{m}{s^2},\unitfrac[3]{m}{s^2}]$. The drag term $c_d$ is assumed to be $\unitfrac[0.3]{m}{s^2}$. The forward velocity is limited according to driving regulations to $x_3 \in [\unitfrac[90]{km}{h},\unitfrac[120]{km}{h}]$ and the lateral velocity $x_4$ is limited to $\pm \unitfrac[20]{km}{h}$. 
Additionally, we impose the lane constraint $x_2 \in [-\unit[3.5]{m},\unit[3.5]{m}]$.
We consider a reference frame moving in longitudinal direction with a constant velocity of $\unitfrac[100]{km}{h}$. All quantities are with respect to this frame. In \reffig{fig:cs:track} the origin of this frame is always $(0,0)$.
A discrete-time version of \eqref{eqn:doubleint} with sampling time $T_s = 0.2$~seconds is used.

A truck is approaching the controlled car from behind with an uncertain initial distance from the controlled car in the range $[\unit[7.5]{m},\unit[18]{m}]$ and an uncertain velocity between $\unitfrac[110]{km}{h}$ and $\unitfrac[120]{km}{h}$. The goal is to escape the approaching truck by performing a lane change. This needs to be accomplished without crashing for any realization of the stated position and velocity uncertainties of the other vehicles.
To model this scenario as an LTL specification, we introduce the three atomic propositions $p_{\rm car 1}$, $p_{\rm car 2}$ and $p_{\rm truck}$ for the obstacles and a proposition $p_{\rm goal}$ for the goal set that we want to reach before being hit by the truck. Hence, we want to satisfy the specification $\varphi := \always\big( \lnot p_{\rm car 1} \land \lnot p_{\rm car 2} \land \lnot p_{\rm truck}\big) \land \eventually \always p_{\rm goal}$ for all realizations $\mathbf{w} \in \setbf{W}$ of the uncertainty. The sets corresponding to the atomic propositions are
\begin{align*}
	\set{P}_{\rm car 1} &:= \big\{ (x,w) \sep{\big} 0 \leq x_1 \leq 6.75 + w_1\,, 0 \leq x_2 \leq 3.5 \big\}\,,\\
	\set{P}_{\rm car 2} &:= \big\{ (x,w) \sep{\big} 31.25 - w_2 \leq x_1 \leq 38\,, 0 \leq x_2 \leq 3.5 \big\}\,,\\
	\set{P}_{\rm truck} &:= \big\{ (x,w) \sep{\big} 0 \leq x_1 \leq w_3\,, -3.5 \leq x_2 \leq 0 \big\}\,,\\
	\set{P}_{\rm goal} &:= \big\{ (x,w) \sep{\big} 0 \leq x_2 \leq 3.5 \big\}\,,
\end{align*}
and the set of disturbances $\setbf{W}$ is defined as
\begin{align*}
	\setbf{W} := \big\{ \mathbf{w} \in \reals{3(N+1)} \sep{\big} & w_{k,i} = \sum_{j=0}^k\nolimits \omega_{j,i} \text{ for } i=1,\ldots,3 \text{, with }\\ 
	&\: \omega_{j,1} \in [-d_c,d_c]\,,\: \omega_{j,2} \in [-d_c,d_c] \text{ for } j=0,\ldots,N\\
	&\: \omega_{0,3} \in [\unit[-13.5]{m},\unit[0]{m}]\,,\\
	&\: \omega_{j,3} \in [T_s \cdot \unitfrac[10]{km}{h},T_s \cdot \unitfrac[20]{km}{h}]  \text{ for } j=1,\ldots,N\: \big\}\,,
\end{align*}
with $d_c := T_s \cdot \unitfrac[6.6]{km}{h}$.
Because a point model is used, the position constraints of all vehicles have to be modified to take into account the car's shape which is $\unit[4.5]{m}$ long and $\unit[2]{m}$ wide. For simplicity, these margins are omitted in the illustration.

The planning horizon is $N=20$ and we use an objective function that promotes \emph{minimum time} solutions, additionally penalizing the control effort.
\begin{equation*}
	J(\mathbf{H},\mathbf{h},\delta) := \underbrace{\sum_{k=0}^N k \delta_{{\rm goal},k}}_{\text{minimum time}} + \gamma\big( \underbrace{\|\mathbf{h}\|_2^2 + \|\mathbf{H}\|_F^2 \max_{\mathbf{w}\in\setbf{W}} \|\mathbf{w}\|_2^2}_{\text{control effort}} \big)\,,
\end{equation*}
with $\gamma = 0.001$. The binary variable $\delta_{{\rm goal},k}$ equals one if the goal is reached at time $k$. We denote by $\|\cdot\|_2$ the 2-norm and by $\|\cdot\|_F$ the Frobenius norm.

We generated both an open-loop and a disturbance feedback policy using YALMIP \cite{lofberg2004} and Gurobi \cite{gurobi} to solve the resulting MIQP on an Intel i7 CPU at 2.8GHz. The open-loop policy problem has $4198$ continuous and $273$ binary variables, and $7401$ constraints. An optimal solution was found in $\unit[0.4]{s}$. For the disturbance feedback policy, the proposed inner approximation was used leading to a problem with $35338$ continuous and $273$ binary variables, and $52401$ constraints. An optimal solution was found in $\unit[17.3]{s}$.

\begin{figure}
	\centering
	\begin{minipage}[b]{0.49\linewidth}
	\centering
	\begin{subfigure}[b]{\linewidth}
		\centering
		\includegraphics[width=\linewidth]{run_6_final.tikz}
	\end{subfigure}\\
	\begin{subfigure}[b]{\linewidth}
		\centering
		\includegraphics[width=\linewidth]{run_13_final.tikz}
	\end{subfigure}\\
	\begin{subfigure}[b]{\linewidth}
		\centering
		\includegraphics[width=\linewidth]{run_20_final.tikz}
	\end{subfigure}
	\caption{Trajectory for time steps $k=\{6,13,20\}$ with the feedback policy applied to the worst-case disturbance realization. The \emph{effective goal} is the feasible part of the goal, additionally taking into account the shape of the controlled car.} \label{fig:cs:track}
	\end{minipage}
	\hspace*{0.005\linewidth}
	\begin{minipage}[b]{0.245\linewidth}
		\centering
		\includegraphics[width=\linewidth]{run_comparison.tikz}
		\caption{Trajectories of controlled car corresponding to random disturbances ({\color{lightred}\textbf{solid}}), best-case ({\color{blue}\textbf{dotted}}) and worst-case ({\color{blue}\textbf{dashed}}) disturbance.} \label{fig:cs:traj}
	\end{minipage}
	\hspace*{0.005\linewidth}
	\begin{minipage}[b]{0.23\linewidth}
		\centering
		\includegraphics[width=\linewidth]{controleffort.tikz}
		\caption{Control effort $\|\mathbf{u}\|_2$ for different disturbance realizations ({\color{blue}\textbf{solid}}), best ({\color{red}\textbf{circle}}), worst ({\color{red}\textbf{square}}) and open-loop case ({\color{blue}\textbf{dashed}}).} \label{fig:cs:obj}
	\end{minipage}
\end{figure}

In \reffig{fig:cs:track} the trajectory resulting from applying the synthesized feedback policy to the worst-case disturbance realization is illustrated for the time steps $k=6,13$ and $20$. The trajectory is feasible. Furthermore, the goal set is reached after $13$ time steps and the car remains there for the remainder of the planning horizon. 
\reffig{fig:cs:traj} illustrates trajectories for different disturbances taken uniformly randomly from $\setbf{W}$ ({\color{lightred}\textbf{solid}}). Additionally, the trajectories for the best- and worst-case disturbance are illustrated ({\color{blue}\textbf{dotted}} and {\color{blue}\textbf{dashed}}, respectively). The open-loop trajectory does not differ substantially from the feedback trajectory resulting from the worst-case disturbance and is therefore omitted. 
Finally in \reffig{fig:cs:obj} we give the overall control effort $\|\mathbf{u}\|_2$ (sorted) that was needed for $1000$ different disturbance realizations ({\color{blue}\textbf{solid line}}), as well as for the best-case ({\color{red}\textbf{circle}}), worst-case ({\color{red}\textbf{square}}) and the open-loop case ({\color{blue}\textbf{dashed line}}). This illustrates that the feedback policy allows the reduction of the control effort compared to the open-loop policy. Furthermore, for all disturbance samples, the feedback policy required $13$ time steps to reach the goal set whereas the open-loop policy required $14$ time steps.

Note that, in practice the controlled car would estimate the disturbances based on measured positions of the vehicles and apply disturbance feedback using these estimates.

\section{Conclusions}
We have addressed the problem of control synthesis for linear systems given bounded LTL specifications affected by uncertain disturbances. We have formulated a robust policy synthesis problem with affine disturbance feedback to synthesize policies that satisfy the specification for all considered disturbances.
We cast this problem as a robust mixed-integer program. Then we introduced a simple inner approximation of the constraint set resulting in an MIQP that can be solved using general-purpose mixed-integer solvers. The proposed method therefore enables the generation of control policies that satisfy the specification robustly and incorporates performance criteria such as minimum time or minimum control effort. The framework was applied to a numerical case study of guiding a car in a lane changing maneuver on a highway.

\bibliographystyle{ieeetr}
\bibliography{references}

\end{document}